\documentclass{amsart}

\usepackage{amsfonts,amssymb,verbatim,amsmath,amsthm,latexsym,textcomp,amscd}
\usepackage{latexsym,amsfonts,amssymb,epsfig,verbatim}
\usepackage{amsmath,amsthm,amssymb,latexsym,graphics,textcomp}
\usepackage{graphicx}
\usepackage{color}
\usepackage{url}

\input xy
\xyoption{all}

\begin{document}

\newtheorem{theorem}{Theorem}[section]
\newtheorem{prop}[theorem]{Proposition}
\newtheorem{lemma}[theorem]{Lemma}
\newtheorem{cor}[theorem]{Corollary}
\newtheorem{definition}[theorem]{Definition}
\newtheorem{conj}[theorem]{Conjecture}
\newtheorem{rmk}[theorem]{Remark}
\newtheorem{claim}[theorem]{Claim}
\newtheorem{defth}[theorem]{Definition-Theorem}

\newcommand{\boundary}{\partial}
\newcommand{\C}{{\mathbb C}}
\newcommand{\integers}{{\mathbb Z}}
\newcommand{\natls}{{\mathbb N}}
\newcommand{\ratls}{{\mathbb Q}}
\newcommand{\bbR}{{\mathbb R}}
\newcommand{\proj}{{\mathbb P}}
\newcommand{\lhp}{{\mathbb L}}
\newcommand{\tube}{{\mathbb T}}
\newcommand{\cusp}{{\mathbb P}}
\newcommand\AAA{{\mathcal A}}
\newcommand\BB{{\mathcal B}}
\newcommand\CC{{\mathcal C}}
\newcommand\DD{{\mathcal D}}
\newcommand\EE{{\mathcal E}}
\newcommand\FF{{\mathcal F}}
\newcommand\GG{{\mathcal G}}
\newcommand\HH{{\mathcal H}}
\newcommand\II{{\mathcal I}}
\newcommand\JJ{{\mathcal J}}
\newcommand\KK{{\mathcal K}}
\newcommand\LL{{\mathcal L}}
\newcommand\MM{{\mathcal M}}
\newcommand\NN{{\mathcal N}}
\newcommand\OO{{\mathcal O}}
\newcommand\PP{{\mathcal P}}
\newcommand\QQ{{\mathcal Q}}
\newcommand\RR{{\mathcal R}}
\newcommand\SSS{{\mathcal S}}
\newcommand\TT{{\mathcal T}}
\newcommand\UU{{\mathcal U}}
\newcommand\VV{{\mathcal V}}
\newcommand\WW{{\mathcal W}}
\newcommand\XX{{\mathcal X}}
\newcommand\YY{{\mathcal Y}}
\newcommand\ZZ{{\mathcal Z}}
\newcommand\CH{{\CC\HH}}
\newcommand\PEY{{\PP\EE\YY}}
\newcommand\MF{{\MM\FF}}
\newcommand\RCT{{{\mathcal R}_{CT}}}
\newcommand\PMF{{\PP\kern-2pt\MM\FF}}
\newcommand\FL{{\FF\LL}}
\newcommand\PML{{\PP\kern-2pt\MM\LL}}
\newcommand\GL{{\GG\LL}}
\newcommand\Pol{{\mathcal P}}
\newcommand\half{{\textstyle{\frac12}}}
\newcommand\Half{{\frac12}}
\newcommand\Mod{\operatorname{Mod}}
\newcommand\Area{\operatorname{Area}}
\newcommand\ep{\epsilon}
\newcommand\hhat{\widehat}
\newcommand\Proj{{\mathbf P}}
\newcommand\U{{\mathbf U}}
 \newcommand\Hyp{{\mathbf H}}
\newcommand\D{{\mathbf D}}
\newcommand\Z{{\mathbb Z}}
\newcommand\R{{\mathbb R}}
\newcommand\Q{{\mathbb Q}}
\newcommand\E{{\mathbb E}}
\newcommand\til{\widetilde}
\newcommand\length{\operatorname{length}}
\newcommand\tr{\operatorname{tr}}
\newcommand\gesim{\succ}
\newcommand\lesim{\prec}
\newcommand\simle{\lesim}
\newcommand\simge{\gesim}
\newcommand{\simmult}{\asymp}
\newcommand{\simadd}{\mathrel{\overset{\text{\tiny $+$}}{\sim}}}
\newcommand{\ssm}{\setminus}
\newcommand{\diam}{\operatorname{diam}}
\newcommand{\pair}[1]{\langle #1\rangle}
\newcommand{\T}{{\mathbf T}}
\newcommand{\inj}{\operatorname{inj}}
\newcommand{\pleat}{\operatorname{\mathbf{pleat}}}
\newcommand{\short}{\operatorname{\mathbf{short}}}
\newcommand{\vertices}{\operatorname{vert}}
\newcommand{\collar}{\operatorname{\mathbf{collar}}}
\newcommand{\bcollar}{\operatorname{\overline{\mathbf{collar}}}}
\newcommand{\I}{{\mathbf I}}
\newcommand{\tprec}{\prec_t}
\newcommand{\fprec}{\prec_f}
\newcommand{\bprec}{\prec_b}
\newcommand{\pprec}{\prec_p}
\newcommand{\ppreceq}{\preceq_p}
\newcommand{\sprec}{\prec_s}
\newcommand{\cpreceq}{\preceq_c}
\newcommand{\cprec}{\prec_c}
\newcommand{\topprec}{\prec_{\rm top}}
\newcommand{\Topprec}{\prec_{\rm TOP}}
\newcommand{\fsub}{\mathrel{\scriptstyle\searrow}}
\newcommand{\bsub}{\mathrel{\scriptstyle\swarrow}}
\newcommand{\fsubd}{\mathrel{{\scriptstyle\searrow}\kern-1ex^d\kern0.5ex}}
\newcommand{\bsubd}{\mathrel{{\scriptstyle\swarrow}\kern-1.6ex^d\kern0.8ex}}
\newcommand{\fsubeq}{\mathrel{\raise-.7ex\hbox{$\overset{\searrow}{=}$}}}
\newcommand{\bsubeq}{\mathrel{\raise-.7ex\hbox{$\overset{\swarrow}{=}$}}}
\newcommand{\tw}{\operatorname{tw}}
\newcommand{\base}{\operatorname{base}}
\newcommand{\trans}{\operatorname{trans}}
\newcommand{\rest}{|_}
\newcommand{\bbar}{\overline}
\newcommand{\UML}{\operatorname{\UU\MM\LL}}
\newcommand{\EL}{\mathcal{EL}}
\newcommand{\tsum}{\sideset{}{'}\sum}
\newcommand{\tsh}[1]{\left\{\kern-.9ex\left\{#1\right\}\kern-.9ex\right\}}
\newcommand{\Tsh}[2]{\tsh{#2}_{#1}}
\newcommand{\qeq}{\mathrel{\approx}}
\newcommand{\Qeq}[1]{\mathrel{\approx_{#1}}}
\newcommand{\qle}{\lesssim}
\newcommand{\Qle}[1]{\mathrel{\lesssim_{#1}}}
\newcommand{\simp}{\operatorname{simp}}
\newcommand{\vsucc}{\operatorname{succ}}
\newcommand{\vpred}{\operatorname{pred}}
\newcommand\fhalf[1]{\overrightarrow {#1}}
\newcommand\bhalf[1]{\overleftarrow {#1}}
\newcommand\sleft{_{\text{left}}}
\newcommand\sright{_{\text{right}}}
\newcommand\sbtop{_{\text{top}}}
\newcommand\sbot{_{\text{bot}}}
\newcommand\sll{_{\mathbf l}}
\newcommand\srr{_{\mathbf r}}
\newcommand\geod{\operatorname{\mathbf g}}
\newcommand\mtorus[1]{\boundary U(#1)}
\newcommand\A{\mathbf A}
\newcommand\Aleft[1]{\A\sleft(#1)}
\newcommand\Aright[1]{\A\sright(#1)}
\newcommand\Atop[1]{\A\sbtop(#1)}
\newcommand\Abot[1]{\A\sbot(#1)}
\newcommand\boundvert{{\boundary_{||}}}
\newcommand\storus[1]{U(#1)}
\newcommand\Momega{\omega_M}
\newcommand\nomega{\omega_\nu}
\newcommand\twist{\operatorname{tw}}
\newcommand\modl{M_\nu}
\newcommand\MT{{\mathbb T}}
\newcommand\Teich{{\mathcal T}}
\renewcommand{\Re}{\operatorname{Re}}
\renewcommand{\Im}{\operatorname{Im}}

\title{Moebius rigidity for simply connected, negatively curved surfaces}

\author{Kingshook Biswas}
\address{Indian Statistical Institute, Kolkata, India. Email: kingshook@isical.ac.in}

\begin{abstract} Let $X, Y$ be complete, simply connected Riemannian surfaces with pinched negative curvature $-b^2 \leq K \leq -1$.
We show that if $f : \partial X \to \partial Y$ is a Moebius homeomorphism between the boundaries at infinity of $X, Y$, then $f$ extends to
an isometry $F : X \to Y$. This can be viewed as a generalization of Otal's marked length spectrum rigidity theorem for closed, negatively curved surfaces,
in the sense that Otal's theorem asserts that if $X, Y$ admit properly discontinuous, cocompact, free actions by groups of isometries and the boundary map $f$ is
Moebius and equivariant with respect to these actions then it extends to an isometry. In our case there are no cocompactness or equivariance assumptions, indeed the isometry
groups of $X, Y$ may be trivial.
\end{abstract}

\bigskip

\maketitle

\tableofcontents

\section{Introduction}

\medskip

We continue in this article the study of Moebius maps between boundaries of CAT(-1) spaces undertaken in \cite{biswas3}, \cite{biswas4}, \cite{biswas5},
\cite{biswas6}, \cite{biswas7}. The principal question is whether a Moebius homeomorphism between the boundaries at infinity of two CAT(-1) spaces
extends to an isometry between the spaces. We recall that the boundary $\partial X$ of a CAT(-1) space comes equipped with a positive function on
the set of quadruples of distinct points in $\partial X$, called the cross-ratio, and a map $f : \partial X \to \partial Y$ between boundaries is
said to be Moebius if it preserves cross-ratios.

\medskip

Bourdon \cite{bourdon2} showed that if $X$ is a rank one symmetric space of noncompact type with the metric normalized so that the maximum of the
sectional curvatures equals $-1$, and $Y$ is any CAT(-1) space, then any Moebius embedding $f : \partial X \to \partial Y$ extends to an
isometric embedding $F : X \to Y$. In \cite{biswas3} it was shown that if $X, Y$ are proper, geodesically complete CAT(-1) spaces, then any Moebius
homeomorphism $f : \partial X \to \partial Y$ extends to a $(1, \log 2)$-quasi-isometry $F : X \to Y$. This extension was shown in \cite{biswas5}
to coincide with a certain geometrically defined extension of Moebius maps called the {\it circumcenter extension}. For $X, Y$ complete,
simply connected Riemannian manifolds of pinched negative curvature $-b^2 \leq K \leq -1$, the main result of \cite{biswas3} was improved in
\cite{biswas5} to show that the circumcenter extension $F : X \to Y$ of a Moebius homeomorphism $f : \partial X \to \partial Y$ is a $(1, (1 - 1/b)\log 2)$-quasi-isometry.
The case of complete, simply connected Riemannian manifolds $X, Y$ of pinched negative curvature $-b^2 \leq K \leq -1$ was further studied in
\cite{biswas6}, where it was shown that if $f : \partial X \to \partial Y$ and $g : \partial Y \to \partial X$ are mutually inverse Moebius homeomorphisms,
then their circumcenter extensions $F : X \to Y$ and $G : Y \to X$ are $\sqrt{b}$-bi-Lipschitz homeomorphisms which are inverses of each other.
Another case which has been considered is that of compact deformations of a negatively curved manifold \cite{biswas4}, \cite{biswas7}. Here, we consider
a complete, simply connected Riemannian manifold $(X, g_0)$ of pinched negative curvature $-b^2 \leq K_{g_0} \leq -1$, and a Riemannian metric $g_1$ on $X$
such that $g_1 = g_0$ outside a compact in $X$, and such that $g_1$ has sectional curvature bounded above by $-1$. The identity map $id : (X, g_0) \to (X, g_1)$ is
bi-Lipschitz, and thus induces a homeomorphism $f : \partial_{g_0} X \to \partial_{g_1} X$ between the boundaries at infinity of $(X, g_0)$ and $(X, g_1)$.
While some partial results were proved in \cite{biswas4}, in \cite{biswas7} a complete solution to the problem in this case was obtained: if the boundary map
$f : \partial_{g_0} X \to \partial_{g_1} X$ is Moebius, then its circumcenter extension $F : (X, g_0) \to (X, g_1)$ is an isometry.

\medskip

In the present article we obtain a complete solution to the problem of extending Moebius maps to isometries
for the case of complete, simply connected Riemannian manifolds
of pinched negative curvature in dimension two:

\medskip

\begin{theorem} \label{mainthm} Let $X, Y$ be complete, simply connected Riemannian surfaces of pinched negative curvature $-b^2 \leq K \leq -1$. If
$f : \partial X \to \partial Y$ is a Moebius homeomorphism, then the circumcenter extension of $f$ is an isometry $F : X \to Y$.
\end{theorem}

\medskip

The above theorem may be viewed as a generalization of the well-known result of Otal on marked length spectrum rigidity for closed, negatively curved surfaces
\cite{otal2}. This result states that if two closed, negatively curved surfaces have the same marked length spectrum, then they are isometric. It is well-known
that two closed, negatively curved manifolds have the same marked length spectrum if and only if there is an equivariant Moebius map between the boundaries
of their universal covers (see \cite{otal1} and section 5 of \cite{biswas3}). Thus Otal's result is equivalent to the following: if $X, Y$ are complete,
simply connected Riemannian surfaces with curvature bounded above by $-1$, admitting free, properly discontinuous, cocompact, isometric actions by a discrete group $\Gamma$,
and $f : \partial X \to \partial Y$ is an equivariant Moebius map, then $f$ extends to an isometry $F : X \to Y$. We remark that the cocompactness of the actions is
crucial to Otal's proof, where a certain invariant is defined by integrating over the compact quotient $X/\Gamma$. In our case, we do not assume existence of any isometric
group actions or equivariance of the Moebius map, indeed the isometry groups of $X, Y$ may well be trivial.

\medskip

The proof of Theorem \ref{mainthm} relies on certain properties of the circumcenter extension proved in \cite{biswas7}. In section 2 we recall the necessary
preliminaries on Moebius maps and the circumcenter extension, and then in section 3 we prove the main theorem.

\medskip

\section{Preliminaries}

\medskip

For details and proofs of the assertions made in this section we refer to \cite{biswas3}, \cite{biswas5}, \cite{biswas6}, \cite{biswas7}.

\medskip

\subsection{Moebius metrics and visual metrics}

\medskip

Let $(Z, \rho_0)$ be a compact metric space of diameter one. For a metric $\rho$ on $Z$, the cross-ratio with respect to the metric $\rho$ is
the function of quadruples of distinct points in $Z$ defined by
$$
[\xi, \xi', \eta, \eta']_{\rho} := \frac{\rho(\xi, \eta)\rho(\xi', \eta')}{\rho(\xi, \eta')\rho(\xi', \eta)}
$$
A metric $\rho$ on $Z$ is said to be antipodal if it has diameter one and for any $\xi \in Z$ there exists $\eta \in Z$ such that $\rho(\xi, \eta) = 1$.
We assume that the metric $\rho_0$ is antipodal. We say that two metrics $\rho_1, \rho_2$ on $Z$ are Moebius equivalent if for all quadruples
of distinct points in $Z$, the cross-ratios with respect to the two metrics are equal. We let $\mathcal{M}(Z, \rho_0)$ denote the set of all antipodal
metrics on $Z$ which are Moebius equivalent to $\rho_0$. For any $\rho_1, \rho_2 \in \mathcal{M}(Z, \rho_0)$, there exists a positive continuous function
on $Z$ called the derivative of $\rho_2$ with respect to $\rho_1$, denoted by $\frac{d\rho_2}{d\rho_1}$, such that
$$
\rho_2(\xi, \eta)^2 = \frac{d\rho_2}{d\rho_1}(\xi)\frac{d\rho_2}{d\rho_1}(\eta) \rho_1(\xi,\eta)^2
$$
for all $\xi, \eta \in Z$, and such that
$$
\frac{d\rho_2}{d\rho_1}(\xi) = \lim_{\eta \to \xi} \frac{\rho_2(\xi, \eta)}{\rho_1(\xi, \eta)}
$$
for all non-isolated points $\xi$ of $Z$. Moreover, 
$$
\left( \max_{\xi \in Z} \frac{d\rho_2}{d\rho_1}(\xi) \right) \cdot \left( \min_{\xi \in Z} \frac{d\rho_2}{d\rho_1}(\xi) \right) = 1
$$
The set $\mathcal{M}(Z, \rho_0)$ admits a natural metric defined by
$$
d_{\mathcal{M}}(\rho_1, \rho_2) = \sup_{\xi \in Z} \log \frac{d\rho_2}{d\rho_1}(\xi)
$$
The metric space $(\mathcal{M}(Z, \rho_0), d_{\mathcal{M}})$ is proper.

\medskip

Let $X$ be a proper, geodesically complete CAT(-1) space with boundary at infinity $\partial X$. For any $x \in X$, there is a metric $\rho_x$ on $\partial X$ called the
visual metric based at $x$, defined by
$$
\rho_x(\xi, \eta) = e^{-(\xi|\eta)_x} \ , \xi, \eta \in \partial X,
$$
where $(\xi|\eta)_x$ is the Gromov inner product between the boundary points $\xi, \eta \in \partial X$ with respect to the basepoint $x$, defined by
$$
(\xi|\eta)_x = \lim_{y \to \xi, z \to \eta} \frac{1}{2} (d(x,y) + d(x,z) - d(y,z))
$$
The metric space $(\partial X, \rho_x)$ is compact, of diameter one, and antipodal. Moreover, for any two points $x, y \in X$, the metrics $\rho_x, \rho_y$ are
Moebius equivalent. Thus the metric space $(\mathcal{M}(\partial X, \rho_x), d_{\mathcal{M}})$ is independent of the choice of $x$, and we denote it by simply
$\mathcal{M}(\partial X)$. The map $X \to \mathcal{M}(\partial X), x \mapsto \rho_x$ is an isometric embedding with image $1/2 \log 2$-dense in $\mathcal{M}(\partial X)$.

\medskip

\subsection{The circumcenter extension}

\medskip

Let $X, Y$ be complete, simply connected Riemannian manifolds of pinched negative curvature $-b^2 \leq K \leq -1$, and suppose there is a
Moebius homeomorphism $f : \partial X \to \partial Y$. The Moebius map $f$ induces a homeomorphism between the unit tangent bundles
$\phi : T^1 X \to T^1 Y$ which conjugates the geodesic flows. The map $\phi$ is defined as follows: given $v \in T^1 X$, let $\gamma : \mathbb{R} \to X$ be the unique
bi-infinite geodesic such that $\gamma'(0) = v$, then let $x = \gamma(0), \xi = \gamma(+\infty), \eta = \gamma(-\infty)$. Let $(f(\xi), f(\eta)) \subset Y$
denote the unique (unparametrized) bi-infinite geodesic in $Y$ with endpoints $f(\xi), f(\eta)$. There exists a unique $y \in (f(\xi), f(\eta))$ such that
$$
\frac{d f_* \rho_x}{d\rho_y}(f(\xi)) = 1
$$
Let $\tilde{\gamma} : \mathbb{R} \to Y$ be the unique bi-infinite geodesic such that $\tilde{\gamma}(0) = y, \tilde{\gamma}(+\infty) = f(\xi),
\tilde{\gamma}(-\infty) = f(\eta)$. We then define $\phi(v) = \tilde{\gamma}'(0) \in T^1 Y$.

\medskip

Recall that in the CAT(-1) space $Y$, any bounded subset $B$ has a unique circumcenter $c(B) \in Y$, which is the unique point minimizing the
function $y \in Y \mapsto \sup_{z \in B} d(y,z)$. Let $(\phi_t : T^1 Y \to T^1 Y)_{t \in \mathbb{R}}$ be the geodesic flow of $Y$.
For any $v \in T^1 Y$, let $p(v) = \gamma(+\infty) \in \partial Y$, where $\gamma : \mathbb{R} \to Y$ is the unique geodesic such that
$\gamma'(0) = v$. This defines a continuous map $p : T^1 Y \to \partial Y$. Let $\pi : T^1 Y \to Y$ denote the canonical projection.
In \cite{biswas5}, it is shown that for any compact
subset $K \subset T^1 Y$ such that $p(K) \subset \partial Y$ is not a singleton, the limit of the circumcenters $c(\pi(\phi_t(K)))$ exists
as $t \to +\infty$. The limit is called the asymptotic circumcenter of the compact $K$ and is denoted by $c_{\infty}(K)$.

\medskip

The circumcenter extension of the Moebius map $f$ is the map $F : X \to Y$ defined by
$$
F(x) = c_{\infty}(\phi(T^1_x X))
$$
In \cite{biswas5}, it is shown that the circumcenter extension $F$ is a $(1, (1-1/b)\log 2)$-quasi-isometry, while in
\cite{biswas6} it is proved that the circumcenter extensions of $f$ and $f^{-1}$ are $\sqrt{b}$-bi-Lipschitz
homeomorphisms which are inverses of each other.

\medskip

For $x \in X$ and $\xi \in \partial X$, let $\overrightarrow{x\xi} \in T^1_x X$ denote the tangent vector $\gamma'(0)$, where
$\gamma$ is the unique geodesic with $\gamma(0) = 0, \gamma(+\infty) = \xi$. For $y \in Y, \eta \in \partial Y$, $\overrightarrow{y\eta} \in T^1_y Y$
is similarly defined. Let $\mu$ be a probability measure on $\partial X$. We say that $\mu$ is balanced at $x \in X$ if
$$
\int_{\partial X} < \overrightarrow{x\xi}, v > d\mu(\xi) = 0
$$
for all $v \in T_x X$. The notion of a probability measure on $\partial Y$ being balanced at a point of $Y$ is similarly defined.

\medskip

Let $F : X \to Y$ be the circumcenter extension of the Moebius map $f : \partial X \to \partial Y$. For $x \in X$, define a
function $u_x$ on $\partial X$ by
$$
u_x(\xi) = \log \frac{d\rho_x}{df^*\rho_{F(x)}}(\xi) \ , \xi \in \partial X.
$$
and let $K_x \subset \partial X$ denote the set where the function $u_x$ achieves its maximum.
In \cite{biswas6}, it is shown that
for any $x \in X$, there exists a probability measure $\mu_x$ on $\partial X$ with support contained in $K_x$
such that $\mu_x$ is balanced at $x \in X$ and $f_* \mu_x$ is balanced at $F(x) \in Y$. We will need the following propositions 
from \cite{biswas7}:

\medskip

\begin{prop} \label{rconstant} (\cite{biswas7}) The function $r : X \to \mathbb{R}$ defined by
$$
r(x) = d_{\mathcal{M}}( \rho_x, f^* \rho_{F(x)} ) \ , x \in X
$$
is constant.
\end{prop}

\medskip

\begin{prop} \label{qisom} (\cite{biswas7}) Let $M \geq 0$ denote the constant value of the function $r$. Then the circumcenter 
map $F : X \to Y$ is a $(1, 2M)$-quasi-isometry.
\end{prop}

\medskip

Given $x \in X$, the flip map $T^1_x X \to T^1_x X, v \mapsto -v$ induces an involution $i_x : \partial X \to \partial X$, defined by requiring 
that $\overrightarrow{xi_x(\xi)} = -\overrightarrow{x\xi}$ for all $\xi \in \partial X$. 

\medskip

\begin{prop} \label{maxmin} (\cite{biswas7}) For $x \in X$, the function $u_x$ achieves its maximum at $\xi \in \partial X$ if and only if 
it achieves its minimum at $i_x(\xi) \in \partial X$. 
\end{prop}

\medskip

\begin{prop} \label{dFstar} (\cite{biswas7}) Let $x \in X$ be a point of differentiability of the circumcenter map $F : X \to Y$. Then 
for any $\xi \in K_x$ and any $v \in T_x X$, we have
$$
< dF_x(v), \overrightarrow{F(x)f(\xi)} > = < v, \overrightarrow{x\xi} >
$$
Equivalently,
$$
dF^*_x( \overrightarrow{F(x)f(\xi)} ) = \overrightarrow{x\xi}
$$
for all $\xi \in K_x$.
\end{prop}

\medskip

The following Lemma follows from Propositions \ref{qisom} and \ref{maxmin}:

\medskip

\begin{lemma} \label{antisom} Suppose for some $x \in X$, there exists $\xi \in \partial X$ such that $\xi, i_x(\xi) \in K_x$. Then the 
circumcenter map $F : X \to Y$ is an isometry.
\end{lemma}

\medskip

\noindent{\bf Proof:} It follows from Proposition \ref{maxmin} that the maximum and minimum values of the function $u_x$ are equal. 
On the other hand we know that the maximum and minimum values are negatives of each other. Since the maximum value equals the constant $M$, 
we have $M = - M$ and hence $M = 0$. It follows from Proposition \ref{qisom} that $F : X \to Y$ is an isometry. $\diamond$
 
\medskip
 
\section{Proof of main theorem}

\medskip

Let $X, Y$ be complete, simply connected Riemannian surfaces of pinched negative curvature $-b^2 \leq K \leq -1$, and let 
$f : \partial X \to \partial Y$ be a Moebius homeomorphism with circumcenter extension $F : X \to Y$. All the tools 
are now in hand for the proof of the main theorem:

\medskip

\noindent{\bf Proof of Theorem \ref{mainthm}:} As mentioned in the previous section, for any $x \in X$ there exists a 
probability measure $\mu_x$ on $\partial X$ with support contained in $K_x$ such that $\mu_x$ is balanced at $x \in X$, and 
$f_* \mu_x$ is balanced at $F(x) \in Y$. As shown in \cite{biswas6}, this is equivalent to the fact that the convex hull in 
$T_x X$ of the compact $\{ \overrightarrow{x\xi} : \xi \in K_x \}$ contains the origin of $T_x X$ and the convex hull in 
$T_{F(x)} Y$ of the compact $\{ \overrightarrow{F(x)f(\xi)} : \xi \in K_x \}$ contains the origin of $T_{F(x)} Y$. 
By the classical Caratheodory theorem on convex hulls, since $X$ is of dimension two this implies that there exists $1 \leq k \leq 3$ and distinct points 
$\xi_1, \dots, \xi_k \in K_x$ and $\alpha_1, \dots, \alpha_k > 0$ (all depending on $x$) such that 
$$
\alpha_1 \overrightarrow{x\xi_1} + \dots + \alpha_k \overrightarrow{x\xi_k} = 0
$$
and $\alpha_1 + \dots + \alpha_k = 1$. Since the vectors $\overrightarrow{x\xi_i}$ are non-zero we must have $2 \leq k \leq 3$. 
Now if any two of the vectors $\overrightarrow{x\xi_i}, \overrightarrow{x\xi_j}$ for some $i \neq j$ are linearly dependent, then 
since they are distinct unit norm vectors we must have $\overrightarrow{x\xi_i} = -\overrightarrow{x\xi_j}$, hence $\xi_j = i_x(\xi_i)$. 
Thus $\xi_i, i_x(\xi_i) \in K_x$, and it follows from Lemma \ref{antisom} that $F$ is an isometry and we are done. In particular if $k = 2$ 
then we are done. Thus we may as well assume 
that for any $x \in X$, there exist distinct points $\xi_1, \xi_2, \xi_3 \in K_x$ and $\alpha_1, \alpha_2, \alpha_3 > 0$ (all depending on $x$) such
that any two of the vectors $\overrightarrow{x\xi_i}, \overrightarrow{x\xi_j}$ for $i \neq j$ are linearly independent. 

\medskip

As mentioned in the previous section, the circumcenter extensions of $f$ and $f^{-1}$ are bi-Lipschitz homeomorphisms which are 
inverses of each other. Thus there are sets $A \subset X$ and $B \subset Y$ of full measure (with respect to the Riemannian volume measures) 
such that $F$ is differentiable at all points
of $A$ and $F^{-1}$ is differentiable at all points of $B$. Since $F$ is bi-Lipschitz, the set $F^{-1}(B) \subset X$ has full measure, 
thus so does the set $C := A \cap F^{-1}(B)$. For any point $x$ of $C$, $F$ is differentiable at $x$, $F^{-1}$ is differentiable at $F(x)$, 
and by the Chain Rule the derivatives $dF_x, dF^{-1}_{F(x)}$ are inverses of each other, so $dF_x : T_x X \to T_{F(x)} Y$ is an 
isomorphism for all $x \in C$. 

\medskip

Now let $x \in C$. As remarked earlier, we may assume that there are distinct points $\xi_1, \xi_2, \xi_3 \in K_x$ such that 
\begin{equation} \label{no1}
\alpha_1 \overrightarrow{x\xi_1} + \alpha_2 \overrightarrow{x\xi_2} + \alpha_3 \overrightarrow{x\xi_3} = 0
\end{equation}
 for some $\alpha_1, \alpha_2, \alpha_3 > 0$, 
and such that any two of the vectors $\overrightarrow{x\xi_i}, \overrightarrow{x\xi_j}$ for $i \neq j$ are linearly independent. 
By Proposition \ref{dFstar}, we have
$$
dF^*_x ( \alpha_1 \overrightarrow{F(x)f(\xi_1)} + \alpha_2 \overrightarrow{F(x)f(\xi_2)} + \alpha_3 \overrightarrow{F(x)f(\xi_3)} ) = \alpha_1 \overrightarrow{x\xi_1} + \alpha_2 \overrightarrow{x\xi_2} + \alpha_3 \overrightarrow{x\xi_3} = 0
$$
and hence
\begin{equation} \label{no2}
\alpha_1 \overrightarrow{F(x)f(\xi_1)} + \alpha_2 \overrightarrow{F(x)f(\xi_2)} + \alpha_3 \overrightarrow{F(x)f(\xi_3)} = 0
\end{equation}
since $dF^*_x$ is an isomorphism. 

\medskip

Now let 
\begin{align*}
a_1 = < \overrightarrow{x\xi_1}, \overrightarrow{x\xi_2} > \ , & \ a_2 = < \overrightarrow{F(x)f(\xi_1)}, \overrightarrow{F(x)f(\xi_2)} > \\ 
b_1 = < \overrightarrow{x\xi_1}, \overrightarrow{x\xi_3} > \ , & \ b_2 = < \overrightarrow{F(x)f(\xi_1)}, \overrightarrow{F(x)f(\xi_3)} > \\
c_1 = < \overrightarrow{x\xi_2}, \overrightarrow{x\xi_3} > \ , & \ c_2 = < \overrightarrow{F(x)f(\xi_2)}, \overrightarrow{F(x)f(\xi_3)} > \\
\end{align*} 
Taking inner products of the left-hand side of equation (\ref{no1}) above with the vectors $\overrightarrow{x\xi_i}, i = 1,2,3$, and 
taking inner products of the left-hand side of equation (\ref{no2}) above with the vectors $\overrightarrow{F(x)f(\xi_i)}, i = 1,2,3$, we find that the 
vectors $\overrightarrow{u_i} := (a_i, b_i, c_i), i = 1,2$ both satisfy the same linear system of equations
$$
T \overrightarrow{u} = \overrightarrow{w}
$$
where $T$ is the $3 \times 3$ matrix 
$$
T = \begin{pmatrix}
     \alpha_2 & \alpha_3 & 0 \\
     \alpha_1 & 0 & \alpha_3 \\
     0 & \alpha_1 & \alpha_2 
     \end{pmatrix}
$$
and $\overrightarrow{w}$ is the column vector
$$
\overrightarrow{w} = \begin{pmatrix}
                     -\alpha_1 \\
                     -\alpha_2 \\
                     -\alpha_3
                     \end{pmatrix}
$$
A computation gives $\det(T) = -2\alpha_1 \alpha_2 \alpha_3 < 0$, so $T$ is nonsingular, and it follows that $\overrightarrow{u_1} = \overrightarrow{u_2}$. 
We thus have 
$$
< \overrightarrow{x\xi_i}, \overrightarrow{x\xi_j} > = < \overrightarrow{F(x)f(\xi_i)}, \overrightarrow{F(x)f(\xi_j)} >
$$
for all $1 \leq i,j \leq 3$. On the other hand, by Proposition \ref{dFstar}, we have 
$$
< \overrightarrow{x\xi_i}, \overrightarrow{x\xi_j} > = < dF_x(\overrightarrow{x\xi_i}), \overrightarrow{F(x)f(\xi_j)} > 
$$
for all $1 \leq i,j \leq 3$, thus
\begin{equation} \label{no3}
< dF_x(\overrightarrow{x\xi_i}), \overrightarrow{F(x)f(\xi_j)} > = < \overrightarrow{F(x)f(\xi_i)}, \overrightarrow{F(x)f(\xi_j)} >
\end{equation}
for all $1 \leq i,j \leq 3$. Since the dimension of $X$ is two, the span of the vectors $\overrightarrow{x\xi_1}, \overrightarrow{x\xi_2}$ equals 
$T_x X$, thus since $dF^*_x$ is an isomorphism it follows from Proposition \ref{dFstar} that the span of the vectors 
$\overrightarrow{F(x)f(\xi_1)}, \overrightarrow{F(x)f(\xi_2)}$ equals $T_{F(x)} Y$. Fixing $i$ and putting $j = 1,2$ in equation (\ref{no3}) above, 
it follows that
$$
dF_x(\overrightarrow{x\xi_i}) = \overrightarrow{F(x)f(\xi_i)}
$$
for all $1 \leq i \leq 3$. Applying $dF^*_x$ to both sides of the above equation and using Proposition \ref{dFstar} it follows that
$$
dF^*_x dF_x ( \overrightarrow{x\xi_i} ) = \overrightarrow{x\xi_i}
$$
for all $1 \leq i \leq 3$. Since the vectors $\overrightarrow{x\xi_1}, \overrightarrow{x\xi_2}$ span $T_x X$, it follows that 
$dF^*_x dF_x = id$. 

\medskip

Thus $dF_x$ is an isometry for all $x \in C$, in particular $||dF_x|| = 1$ for all $x \in C$. Now it is a classical fact that if a 
Lipschitz map $F$ between Riemannian manifolds satisfies $||dF|| \leq L$ almost everywhere for some constant $L$, then $F$ is $L$-Lipschitz. 
It follows that the circumcenter extension $F : X \to Y$ of the Moebius map $f$ is $1$-Lipschitz. Now we know from \cite{biswas6} that $F^{-1}$ is the circumcenter 
extension of the Moebius map $f^{-1}$. Applying the same argument to the Moebius map $f^{-1}$, we obtain that its circumcenter extension $F^{-1}$ is also 
$1$-Lipschitz. Since both $F$ and $F^{-1}$ are $1$-Lipschitz, $F$ is an isometry. $\diamond$

\medskip

\bibliography{moeb}
\bibliographystyle{alpha}

\end{document}